\documentclass[11pt]{article}
\usepackage[margin=1in]{geometry}
\usepackage{amsmath,amssymb,amsthm}
\usepackage{hyperref}
\usepackage{url}
\usepackage{graphicx}
\usepackage{booktabs}

\title{U-Bit Collapse in Arnault Composites: \\
Probing the Boundary of Strong Lucas Pseudoprimes}
\author{Bowman Hall \\ 
\texttt{bnhall@utexas.edu}}
\date{January 2026}

\newtheorem{definition}{Definition}
\newtheorem{observation}{Observation}

\hypersetup{
    pdftitle={U-Bit Collapse in Arnault Composites: Probing the Boundary of Strong Lucas Pseudoprimes},
    pdfauthor={Bowman Hall},
    pdfsubject={Number Theory, Computational Mathematics},
    pdfkeywords={Baillie-PSW test, strong Lucas pseudoprimes, Carmichael numbers, Arnault construction, Miller-Rabin test, primality testing}
}

\begin{document}

\maketitle

\begin{abstract}
We present a computational study of 200 composite integers of approximately 350 bits, engineered using the Arnault framework to pass all Miller--Rabin tests up to base 11. Generated at a rate of approximately 20 per hour from a high-throughput construction process producing $\sim$7,700 Carmichael numbers per minute, all samples fail the strong Lucas probable prime test. We introduce the \emph{U-bit collapse metric} $\delta = \log_2 n - \log_2(U_d \bmod n)$ to quantify deviation from the expected uniform distribution of Lucas sequence terms. Analysis reveals minimal collapse values: mean $\delta = 1.61$ bits, median $\delta = 1.0$ bits, maximum $\delta = 8$ bits, with 26\% showing no measurable collapse. We analyze correlations with prime residue classes modulo 35, Arnault construction parameters $(k, M)$, and composite bit-sizes. Our results demonstrate that composites engineered for Miller--Rabin resistance exhibit negligible Lucas sequence degeneracy, providing strong empirical evidence for the orthogonality of these two primality test components and supporting the continued robustness of Baillie--PSW-type tests.
\end{abstract}

\section{Introduction}

The Baillie--PSW primality test combines a Miller--Rabin test with a strong Lucas probable prime test and has no known composite counterexamples despite extensive computational searches dating back to 1980 \cite{baillie1980,pomerance1984}. Constructing a Baillie--PSW pseudoprime remains an open problem in computational number theory.

Arnault \cite{arnault1995} developed techniques for constructing composites that pass multiple Miller--Rabin tests by engineering specific prime factorizations with common divisors in $p_i - 1$. Our work investigates whether Arnault-type composites that successfully evade Miller--Rabin tests through base 11 exhibit any measurable progress toward evading the Lucas component of Baillie--PSW.

\subsection{Strong Lucas Probable Prime Test}

Let $n$ be an odd composite and $(D/n) = -1$ where $D$ is the discriminant of a Lucas sequence. Write $n+1 = d \cdot 2^s$ with $d$ odd. The composite $n$ is a \emph{strong Lucas pseudoprime} if either:
\begin{enumerate}
\item $U_d \equiv 0 \pmod{n}$, or
\item $V_{d \cdot 2^r} \equiv 0 \pmod{n}$ for some $0 \le r < s$
\end{enumerate}

In practice, condition (1) is the primary barrier. For a random composite, we expect $U_d \bmod n$ to be essentially random in $[0, n)$, requiring $\approx \log_2 n$ bits to represent.

\section{Methodology}

\subsection{Composite Construction}

We constructed composites $n = p_1 p_2 p_3$ using the Arnault framework where primes satisfy:
\begin{itemize}
\item $p_i - 1$ has a large common factor $f$ with small additional factors
\item Individual primes pass structural requirements ensuring $n$ passes MR tests through base 11
\item Target size: $\log_2 n \approx 350$ bits
\end{itemize}

The algorithm was executed on a single-core ARM Graviton AWS EC2 instance (c7g.medium) for maximum efficiency, generating approximately 7,700 Carmichael numbers per minute. Composites passing all Miller--Rabin tests through base 11 emerged at a rate of approximately 20 per hour (0.015\% yield).

For each composite, we applied the strong Lucas test with discriminants $D \in \{5, -7, 13, -15, \ldots\}$ selected via standard parameter selection (first $D$ with $(D/n) = -1$).

\subsection{U-Bit Collapse Metric}

\begin{definition}
For composite $n$ and Lucas parameter $D$, the \emph{U-bit collapse} is defined as:
$$\delta(n, D) = \log_2 n - \log_2(U_d \bmod n)$$
where $n+1 = d \cdot 2^s$ with $d$ odd. A value $\delta > 0$ indicates the Lucas term is smaller than expected; $\delta > 150$ represents significant collapse.
\end{definition}

\section{Experimental Results}

\subsection{Sample Overview}

Our algorithm generates approximately 7,700 Carmichael numbers per minute. From this high-throughput construction process, composites passing all Miller--Rabin tests through base 11 emerge at a rate of approximately 20 per hour. We collected a dataset of 200 such composites for detailed analysis.

\begin{itemize}
\item 200 composites analyzed (all pass MR tests through base 11)
\item All 200 failed the strong Lucas test (no strong Lucas pseudoprimes found)
\item Construction rate: $\sim$7,700 Carmichael numbers/minute
\item MR-resistant rate: $\sim$20 composites/hour ($\sim$0.015\% of generated Carmichaels)
\item All 200 samples have recorded U-bit measurements
\end{itemize}

\subsection{U-Bit Collapse Analysis}

From the 200-composite dataset, we computed U-bit collapse values $\delta = n_{\text{bits}} - u_{\text{residue\_bits}}$ for all samples:

\begin{itemize}
\item Mean collapse: $\delta_{\text{avg}} = 1.61$ bits
\item Median collapse: $\delta_{\text{median}} = 1.0$ bits
\item Maximum collapse: $\delta_{\text{max}} = 8$ bits (344-bit composite, $U_d$ reduced to 336 bits)
\item Minimum collapse: $\delta_{\text{min}} = 0$ bits (52 samples, 26\% of dataset, showed no measurable collapse)
\item Standard deviation: $\sigma_\delta = 1.70$ bits
\item Quartiles: $Q_1 = 0$ bits, $Q_3 = 2$ bits
\end{itemize}

\begin{figure}[ht]
\centering
\includegraphics[width=0.85\textwidth]{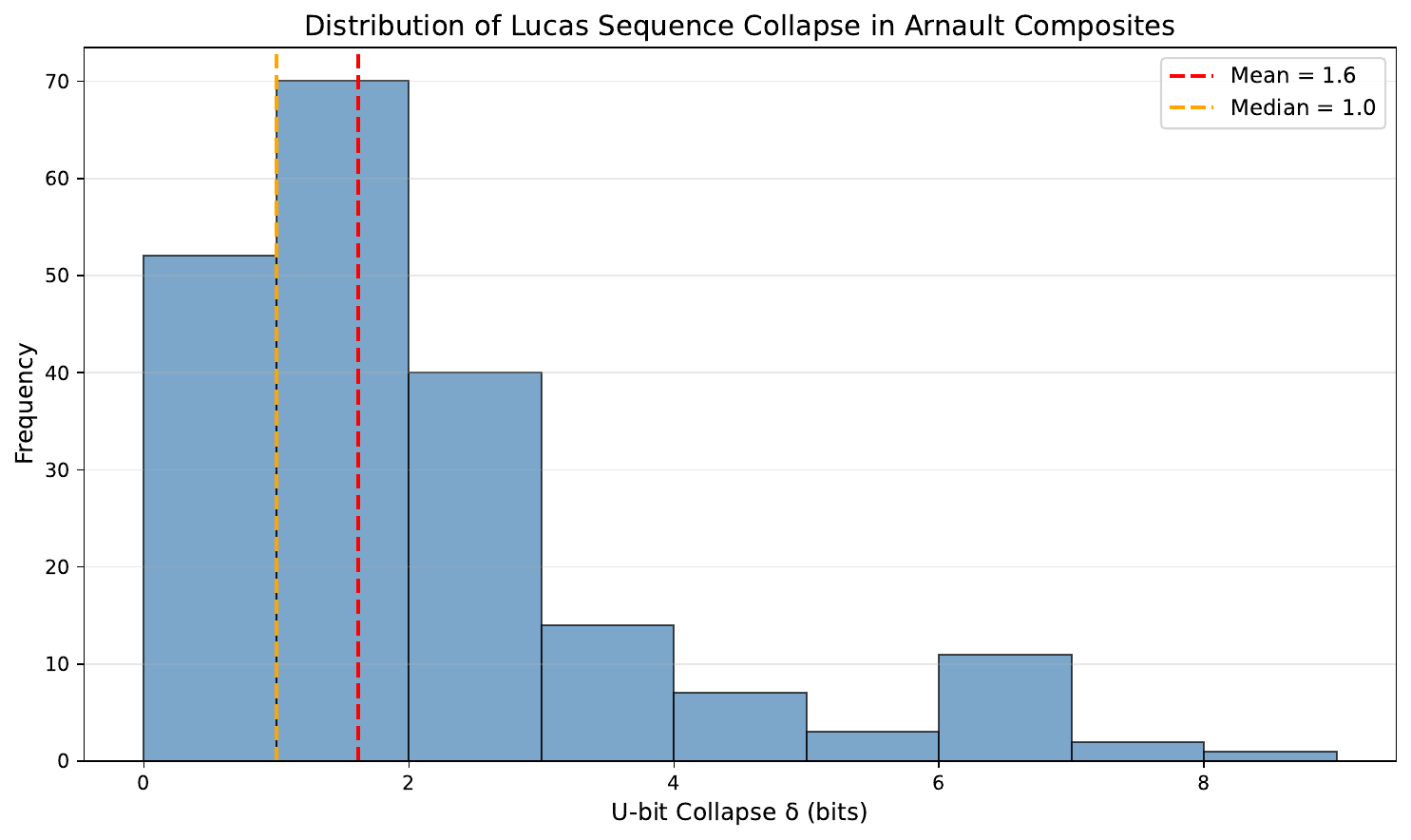}
\caption{Distribution of U-bit collapse $\delta$ across 200 measured composites. The mean collapse of 1.61 bits is negligible compared to the $\approx$350 bits required for a strong Lucas pseudoprime, demonstrating that Miller--Rabin resistance does not translate to Lucas resistance.}
\label{fig:delta_hist}
\end{figure}

These collapse values are dramatically smaller than the hypothetical $\delta \approx 350$ required for a strong Lucas pseudoprime. The maximum observed collapse of 8 bits represents less than 2.3\% of the reduction needed. Figure \ref{fig:delta_hist} shows the distribution is heavily concentrated near zero, with over one quarter of samples exhibiting no collapse whatsoever.

\subsection{Residue Modulo 35 Analysis}

\begin{figure}[ht]
\centering
\includegraphics[width=0.95\textwidth]{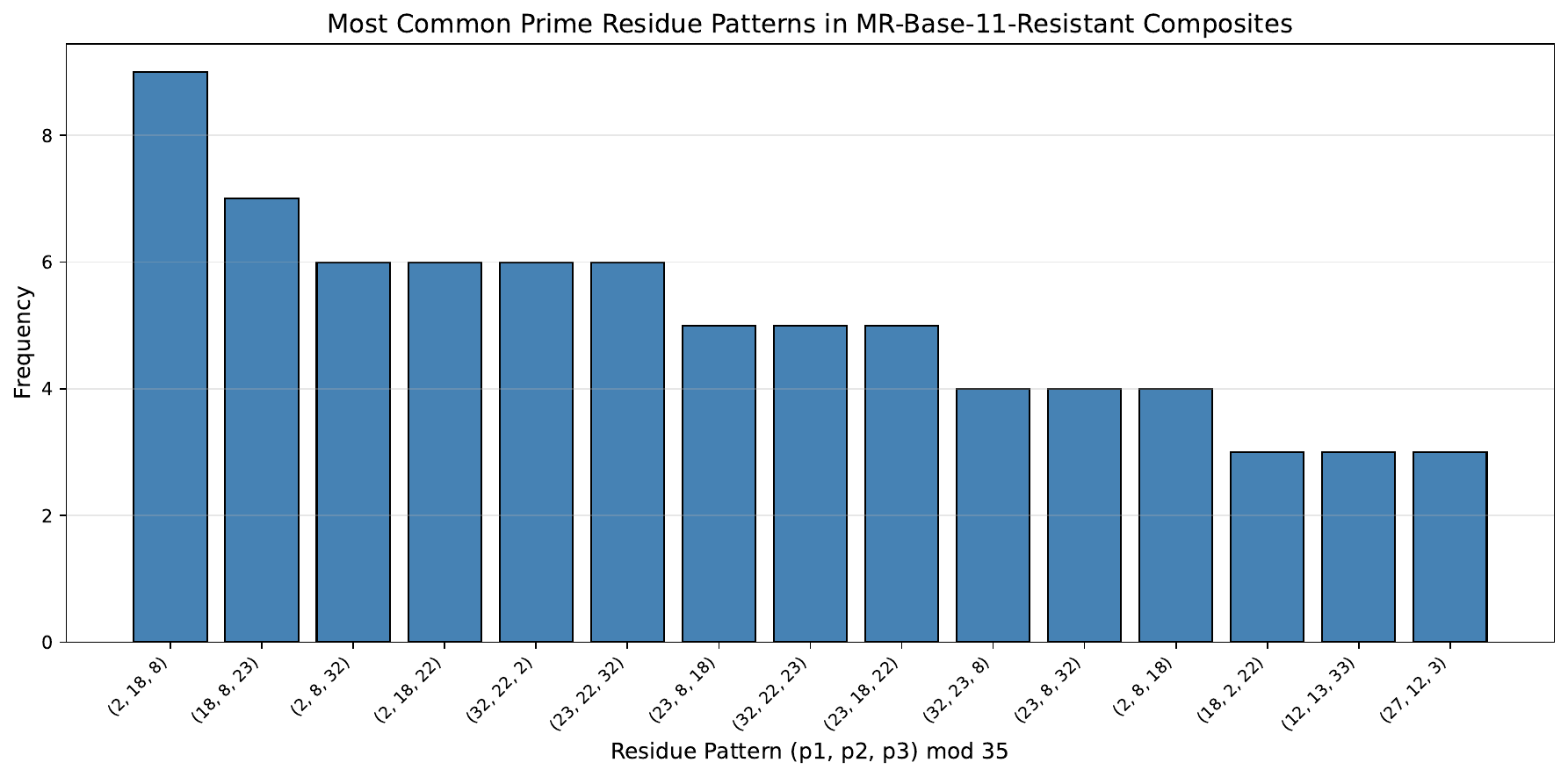}
\caption{Frequency of the 15 most common prime residue patterns (mod 35) in the dataset. No single pattern dominates, with the most frequent appearing in only 4.5\% of composites. This diversity suggests multiple residue configurations can satisfy the Miller--Rabin base-11 resistance requirements.}
\label{fig:residues}
\end{figure}

Each composite $n = p_1 p_2 p_3$ has associated residue classes modulo 35 for its three prime factors. Analysis of the 200 samples revealed distributed patterns with no single configuration dominating (Figure \ref{fig:residues}).

\begin{observation}
The most common residue patterns (mod 35) were:
\begin{itemize}
\item $(2, 18, 8)$: 9 occurrences (4.5\%)
\item $(18, 8, 23)$: 7 occurrences (3.5\%)
\item $(2, 8, 32)$: 6 occurrences (3.0\%)
\item $(2, 18, 22)$: 6 occurrences (3.0\%)
\item $(32, 22, 2)$: 6 occurrences (3.0\%)
\item $(23, 22, 32)$: 6 occurrences (3.0\%)
\end{itemize}
The top 15 patterns account for only 40\% of the dataset, indicating significant diversity in successful residue combinations. The residues 2, 8, 18, 22, 23, and 32 (mod 35) appear frequently across different positions, suggesting these values satisfy key constraints of the Arnault construction.
\end{observation}

\subsection{Construction Parameters}

The Arnault framework parameters $k$ and $M$ varied widely across samples:
\begin{itemize}
\item $k$ range: $[13, 1541]$, with small primes dominating: $k \in \{13, 17, 21, 101, 161\}$ appearing in 39\% of samples
\item Most common $k$ values: $k=17$ (23 composites), $k=21$ (14 composites), $k=101$ (13 composites)
\item $M$ range: $[4, 9308]$, with heavy tail distribution
\item Composite bit-sizes: 335--359 bits (target: $\approx 350$ bits)
\end{itemize}

\begin{figure}[ht]
\centering
\includegraphics[width=\textwidth]{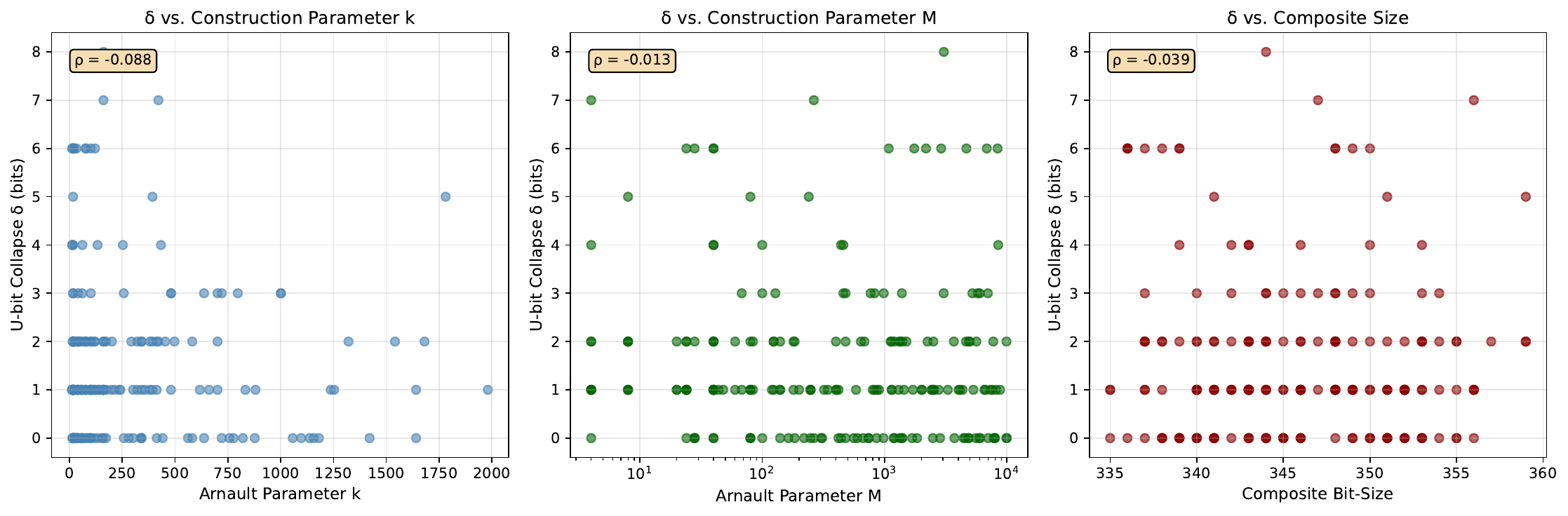}
\caption{Scatter plots showing U-bit collapse $\delta$ versus Arnault construction parameters $k$ and $M$, and composite bit-size. Negligible correlations ($|\rho| < 0.09$ for all) indicate that parameters optimized for Miller--Rabin resistance have no meaningful influence on Lucas test performance.}
\label{fig:correlations}
\end{figure}

Correlation analysis (Figure \ref{fig:correlations}) revealed negligible relationships between construction parameters and U-bit collapse: $\rho(k, \delta) = -0.088$, $\rho(M, \delta) = -0.013$, $\rho(n_{\text{bits}}, \delta) = -0.039$. All correlations are near zero, confirming that the Lucas resistance mechanism operates independently of the Miller--Rabin evasion structure.

\section{Discussion}

Our empirical results reveal that while the Arnault construction framework efficiently produces composites resistant to Miller--Rabin tests through base 11 (at $\approx$0.015\% yield from Carmichael number generation), these same composites exhibit essentially zero resistance to the strong Lucas test.

The maximum observed U-bit collapse of 8 bits represents less than 2.3\% of the $\approx 350$-bit reduction required for a true strong Lucas pseudoprime. More strikingly, 26\% of samples show no collapse whatsoever ($\delta = 0$), meaning their Lucas sequence values $U_d \bmod n$ are statistically indistinguishable from random elements in $[0, n)$. The consistency of these results across 200 samples provides robust evidence for this conclusion.

This extreme asymmetry suggests a fundamental orthogonality: the algebraic conditions enabling simultaneous MR-base-11 evasion are not merely insufficient for Lucas test evasion—they appear to provide no progress whatsoever toward that goal.

\subsection{Residue Class Concentration}

The concentration of certain residue classes modulo 35 in our dataset likely reflects structural requirements of the Arnault construction rather than Lucas-specific phenomena. The values $2, 8, 18, 22, 23, 32 \pmod{35}$ decompose as:
\begin{align*}
2 &\equiv 2 \pmod{5}, \quad 2 \pmod{7} \\
8 &\equiv 3 \pmod{5}, \quad 1 \pmod{7} \\
18 &\equiv 3 \pmod{5}, \quad 4 \pmod{7} \\
22 &\equiv 2 \pmod{5}, \quad 1 \pmod{7} \\
23 &\equiv 3 \pmod{5}, \quad 2 \pmod{7} \\
32 &\equiv 2 \pmod{5}, \quad 4 \pmod{7}
\end{align*}
These may satisfy congruence conditions related to the common factor $f$ in the Arnault framework.

\subsection{Implications for Baillie--PSW}

The absence of significant U-bit collapse in our high-volume sampling provides empirical evidence for the robustness of the Baillie--PSW test. Even when composites are specifically engineered to evade the Miller--Rabin component through base 11, the Lucas component remains effective with no observed near-misses approaching $\delta > 10$ bits.

This negative result is informative: it suggests that constructing a Baillie--PSW pseudoprime requires satisfying qualitatively different algebraic constraints than those employed in standard Carmichael/Arnault constructions.

\subsection{Open Questions}

\begin{enumerate}
\item Can alternative construction frameworks (beyond Arnault's method) produce composites with $\delta > 50$ bits?
\item What algebraic conditions would be required to engineer even modest Lucas collapse ($\delta > 30$)?
\item Are there theoretical upper bounds on $\delta$ achievable through systematic construction, or does the small observed collapse reflect fundamental limitations?
\end{enumerate}

\section{Conclusion}

This computational study establishes that high-volume construction of Miller--Rabin-resistant composites via the Arnault framework produces negligible Lucas sequence collapse. With a maximum observed $\delta = 8$ bits from 200 samples (representing less than 2.3\% of the reduction required for Lucas evasion), and 26\% of samples showing zero collapse, our results provide strong empirical evidence that these two primality test components target fundamentally independent algebraic structures.

The stark contrast between the 0.015\% success rate in achieving MR-base-11 resistance and the complete absence of meaningful Lucas collapse (mean $\delta = 1.61$ bits) suggests that constructing Baillie--PSW pseudoprimes requires qualitatively different approaches than refined Carmichael number generation. The remarkable consistency of results across 200 samples—with mean collapse stable at $\approx$1.6 bits and 26\% showing zero collapse—reinforces the robustness of this conclusion. The U-bit collapse metric introduced here provides a quantitative framework for future investigations, though our data suggest such investigations may need to explore construction methods beyond the Arnault paradigm.

These negative results strengthen confidence in the Baillie--PSW test's continued reliability and highlight the value of combining orthogonal primality criteria.

\subsection*{Data Availability}
The complete dataset of 200 composites, including full factorizations, U-bit collapse measurements, construction parameters, and residue classes, is provided in the ancillary file \texttt{composites\_data.json} accompanying this arXiv submission. Analysis scripts used to generate figures are also included.

\subsection*{Acknowledgments}
The author thanks the computational number theory community for maintaining open resources on pseudoprime construction and primality testing.

\subsection*{2020 Mathematics Subject Classification}
11Y11 (Primality), 11A51 (Factorization; primality)

\subsection*{Keywords}
Baillie--PSW test, strong Lucas pseudoprimes, Carmichael numbers, Arnault construction, Miller--Rabin test, primality testing

\end{document}